%% file: elsarticle-main.tex
\documentclass[review]{elsarticle}

\usepackage[margin=2.5cm]{geometry}

\usepackage{lineno,hyperref}
\modulolinenumbers[5]

\usepackage{multirow}
\usepackage{xcolor}
\usepackage{amsmath}
\usepackage{amsfonts}
\usepackage{stmaryrd}
\usepackage{array}
\usepackage{caption}
\newcolumntype{M}[1]{>{\centering\arraybackslash}m{#1}}

\DeclareMathOperator{\diam}{diam}
\usepackage{algorithm}
\usepackage{algorithmic}
\usepackage{dirtytalk}
\usepackage{todonotes}
\usepackage{svg}

\journal{XXX}

\begin{document}

\begin{frontmatter}

\title{Agglomeration of Polygonal Grids using Graph Neural Networks with applications to Multigrid solvers} 


\author[mymainaddress]{P. F. Antonietti\fnref{myfootnote}}
\ead{paola.antonietti@polimi.it}

\author[mymainaddress]{N. Farenga}
\ead{nicola.farenga@mail.polimi.it}

\author[mymainaddress]{E. Manuzzi\fnref{myfootnote}\corref{mycorrespondingauthor}}
\ead{enrico.manuzzi@mail.polimi.it}

\author[mymainaddress]{G. Martinelli}
\ead{gabriele2.martinelli@mail.polimi.it}

\author[mymainaddress]{L. Saverio}
\ead{luca.saverio@mail.polimi.it}

\address[mymainaddress]{MOX, Department of Mathematics, Politecnico di Milano, p.zza Leonardo da Vinci, 32, I-20133 Milano, Italy}
\fntext[myfootnote]{P. F. Antonietti and E. Manuzzi are members of INDAM-GNCS. P.F. Antonietti has been partially supported by the Ministero dell’Università e della Ricerca [PRIN grant numbers 201744KLJL and 20204LN5N5].}

\cortext[mycorrespondingauthor]{Corresponding author}

\tnotetext[mytitlenote]{Abbreviations: Machine Learning (ML), Graph Neural Networks (GNNs), Polygonal Discontinuous Galerkin (PolyDG), Finite Element Methods (FEMs), MultiGrid (MG).}

\begin{abstract}
\input{sections/abstract}
\end{abstract}


\end{frontmatter}


\input{sections/introduction}
\input{sections/agglom_strategies}
\input{sections/GNN}

\input{sections/validation_grids}
\input{sections/multigrid}
\input{sections/conclusions}

\newpage
\subsection*{CRediT authorship contribution statement}
\noindent \textbf{P. F. Antonietti:} Conceptualization, Methodology, Resources, Writing - review and editing, Supervision, Project administration, Funding acquisition.\\
\textbf{N. Farenga, G. Martinelli, L. Saverio:} Methodology, Software,  Validation, Formal analysis, Investigation, Data curation, Visualization.\\
\textbf{E. Manuzzi:} Conceptualization, Methodology, Software, Validation, Formal analysis, Investigation, Data curation,    Visualization, Supervision, Writing – original draft.

\subsection*{Declaration of competing interest}
The authors declare that they have no known competing financial interests or personal relationships that could have appeared to influence the work reported in this paper.

\subsection*{Acknowledgements}
Funding: P.F. Antonietti has been partially supported by the Ministero dell’Università e della Ricerca [PRIN grant numbers 201744KLJL and 20204LN5N5]. P. F. Antonietti and E. Manuzzi are members of INDAM-GNCS.\\

\clearpage
\bibliography{main_bibfile}

\end{document}

%% file: sections/abstract.tex
Agglomeration-based strategies are important both within adaptive refinement algorithms and to construct scalable multilevel algebraic solvers. In order to automatically perform agglomeration of polygonal grids, we propose the use of Machine Learning (ML) strategies, that can naturally exploit geometrical information about the mesh in order to preserve the grid quality, enhancing performance of numerical methods and reducing the overall computational cost. In particular, we employ the k-means clustering algorithm and Graph Neural Networks (GNNs) to partition the connectivity graph of a computational mesh. Moreover, GNNs have high online inference speed and the advantage to process naturally and simultaneously both the graph structure of mesh and the geometrical information, such as the areas of the elements or their barycentric coordinates. These techniques are compared with METIS, a standard algorithm for graph partitioning, which is meant to process only the graph information of the mesh. We demonstrate that performance in terms of quality metrics is enhanced for ML strategies. Such models also show a good degree of generalization when applied to more complex geometries, such as brain MRI scans, and the capability of preserving the quality of the grid. The effectiveness of these strategies is demonstrated also when applied to MultiGrid (MG) solvers in a Polygonal Discontinuous Galerkin (PolyDG) framework. In the considered experiments, GNNs show overall the best performance in terms of inference speed, accuracy and flexibility of the approach.


%% file: sections/introduction.tex
\section{Introduction}
Many applications in the fields of Engineering and Applied Sciences, such as fluid-structure interaction problems, flow in fractured porous media, and crack and wave propagation problems, are characterized by a strong complexity of the physical domain, possibly involving moving geometries, heterogeneous media, immersed interfaces and complex topographies. Whenever classical Finite Element Methods (FEMs) are employed to discretize the underlying differential model, the process of grid generation can be the bottleneck of the whole simulation, as computational meshes can be composed only of tetrahedral, hexahedral, or prismatic elements. To overcome this limitation, in the last years, there has been a great interest in developing FEMs that can employ general polygons and polyhedra as grid elements for the numerical discretizations of partial differential equations. We mention the mimetic finite difference method \cite{hyman1997numerical,brezzi2005family,brezzi2005convergence,da2014mimetic}, the hybridizable discontinuous Galerkin method \cite{cockburn2008superconvergent,cockburn2009superconvergent,cockburn2009unified,cockburn2010projection}, the Polyhedral Discontinuous Galerkin (PolyDG) method \cite{hesthaven2007nodal,bassi2012flexibility,antonietti2013hp,cangiani2014hp,antonietti2016review,cangiani2017hp,antonietti2021high}, the Virtual Element Method (VEM) \cite{beirao2013basic,beirao2014hitchhiker,beirao2016virtual,da2016mixed,beirao2021recent,book.vem.sema.simai.2022} and the Hybrid High-Order method \cite{di2014arbitrary,di2015hybrid,di2015hybrid2,di2016review,di2019hybrid}.
This calls for the need to develop effective algorithms to handle polygonal and polyhedral (polytopal, for short) grids and to assess their quality (see e.g. \cite{attene2019benchmark}).
For a comprehensive overview we refer to the monographs and special issues \cite{da2014mimetic,cangiani2017hp,di2016review,di2021polyhedral,beirao2021recent,book.vem.sema.simai.2022} and the references therein.
Among the open problems, there is the issue of efficiently handling polytopal mesh agglomeration, i.e., merging mesh elements to obtain coarser grids \cite{chan1998agglomeration,antonietti2020agglomeration,bassi2012flexibility,pan2022agglomeration,gilbert1998geometric}. 
Mesh agglomeration can be used to obtain a coarser discretization of the differential problem at hand, in order to reduce the number of degrees of freedom where not needed and therefore also the computational effort. This operation can be naturally performed in the context of polygonal and polyhedral grids, because of the flexibility in the definition of the shape of mesh elements. This approach has multiple applications in the numerical solution of partial differential equations, for example:
\begin{itemize}
    \item it can be used, with adaptive procedures, to reduce the number of degrees of freedom where not needed because in certain parts of the domain the error is already under control;
    \item it can be used to generate a hierarchy of (nested) coarser grids starting from a fine mesh of a complex physical domain of interest, in order to employ them in multigrid solvers \cite{bassi2012flexibility,bassi2012agglomeration,antonietti2015multigrid,antonietti2017multigrid,xu2017algebraic,chan1998multilevel,multigridV} to accelerate the converge of iterative algebraic;
    \item it can be employed together with domain decomposition techniques \cite{antonietti2007schwarz,antonietti2014domain,toselli2004domain,feng2001two} to obtain a meaningful decomposition of the domain, starting from a fine discretization.
\end{itemize}
Grid agglomeration is a topic quite unexplored because it is not possible to develop such kind of strategies within the framework of classical FEMs. During this operation, it is important to preserve the quality of the underlying mesh, since this might affect the overall performance of the method in terms of stability and accuracy. Indeed a suitably adapted mesh may allow achieving the same accuracy with a much smaller number of degrees of freedom when solving the numerical problem, hence saving memory and computational power. However, since in such a general framework mesh elements may have any shape, there are no well-established strategies to achieve effective agglomeration with an automatic, fast and simple approach.\\
\newline
In recent years there has been a great development of Machine Learning (ML) algorithms, a framework which allows extracting information automatically from data, to enhance and accelerate numerical methods for scientific computing \cite{raissi2019physics,raissi2018hidden,regazzoni2019machine,regazzoni2020machine,hesthaven2018non,ray2018artificial,antonietti2021accelerating,regazzoni2021machine,ANTONIETTI2022110900,antonietti2022machine}.
In this work, we propose to use ML-based strategies to efficiently handle polygonal mesh agglomeration, in order to fully exploit all of the benefits of the above mentioned numerical methods, such as geometrical flexibility and convergence properties. ML techniques can process naturally geometric information about mesh elements and, as a result, they can preserve the geometrical quality of the initial grid during agglomeration, enhancing the performance of numerical methods and reducing the memory storage and the  associated computational cost. The core concept lies in learning the "shape" of mesh elements in order to perform the desired operations accordingly. Such learning needs to be performed in an automatic and flexible way, because of the too high variability of geometries of interest, tailoring the approach to a wide range of different possible situations. Rather than simply trying to decide a priori criteria to perform agglomeration, which would inevitably result in poor performance or high computational cost due to the impossibility of capturing all of the possible situations, ML strategies exploit and process automatically the huge amount of available data to learn only the distribution of the features of interest for the application, leading to high performance and computational efficiency. By combining the a priori approach of classical numerical methods, with the a posteriori approach of ML strategies, it is possible not only to boost existing algorithms but also to develop new algorithms capable to work in more general frameworks.
In particular, we propose the use of the k-means clustering algorithm \cite{kmeans,hartigan1979algorithm,likas2003global,bello2012adaptive} and Graph Neural Networks (GNNs) \cite{1706.02216,gatti1,gatti2,xu2021graph,lecun2015deep}, the latter being deep learning architectures specifically meant to work with graph-structured data. The problem of mesh agglomeration can be re-framed as a graph partitioning problem, by exploiting the connectivity structure of the mesh. In particular, the graph representation of the mesh is obtained by assigning a node to each element of the mesh and connecting with an edge the pair of nodes which are relative to adjacent elements in the original mesh. By exploiting such a representation, these ML techniques can be applied to solve a node classification problem, where each element is assigned to a cluster, which corresponds to an element of the agglomerated mesh. In particular, GNNs have high online inference speed and are able to process naturally and simultaneously both the graph structure of mesh and the geometrical information that can be attached to the nodes, such as the elements areas or their barycentric coordinates, while the k-means can process only the geometrical information. These techniques are compared with METIS \cite{metis}, a standard solver for graph partitioning, that is meant to process only the graph information about the mesh.
To investigate the capabilities of the proposed approaches, we consider a second-order model problem discretized by the PolyDG method in a multigrid framework. We measure effectiveness through an analysis of quality metrics and the number of iterations of the algebraic iterative solver. We also consider the generalization capabilities on a mesh coming from a human brain MRI scan section.\\
\newline
The paper is organized as follows. In Section 2 we present possible agglomeration criteria for polytopes. In Section 3 we propose a general framework to perform mesh agglomeration using ML strategies, with a focus on GNNs. In Section 4 we measure the effectiveness of the proposed agglomeration strategies in terms of computational cost, quality metrics and generalization capabilities over unseen complex physical domains. In Section 5 we present some computations obtained by applying the considered agglomeration strategies to multigrid solvers in a PolyDG framework. In Section 6 we draw some conclusions.

%% file: sections/agglom_strategies.tex
\section{Mesh agglomeration strategies}
\label{sec:graphpart}
We recall that the problem of mesh agglomeration can be re-framed as a graph partitioning problem, by exploiting the connectivity structure of the mesh. In particular, the graph representation of the mesh is obtained by assigning a node to each element of the mesh, and connecting with an edge the pair of nodes which are relative to adjacent elements in the original mesh, i.e. polygons that share at least one edge.
Moreover, features can be assigned to each node, storing geometrical information such as the area of the element or its barycentric coordinates. Finally, partitioning the nodes of the mesh into proper clusters using a suitable algorithm allows to obtain an agglomerated representation of the original mesh.

\subsection{Graph Partitioning}

\label{sec: graph partitioning}

\begin{figure}
    \centering
    \includegraphics[width=.5\textwidth]{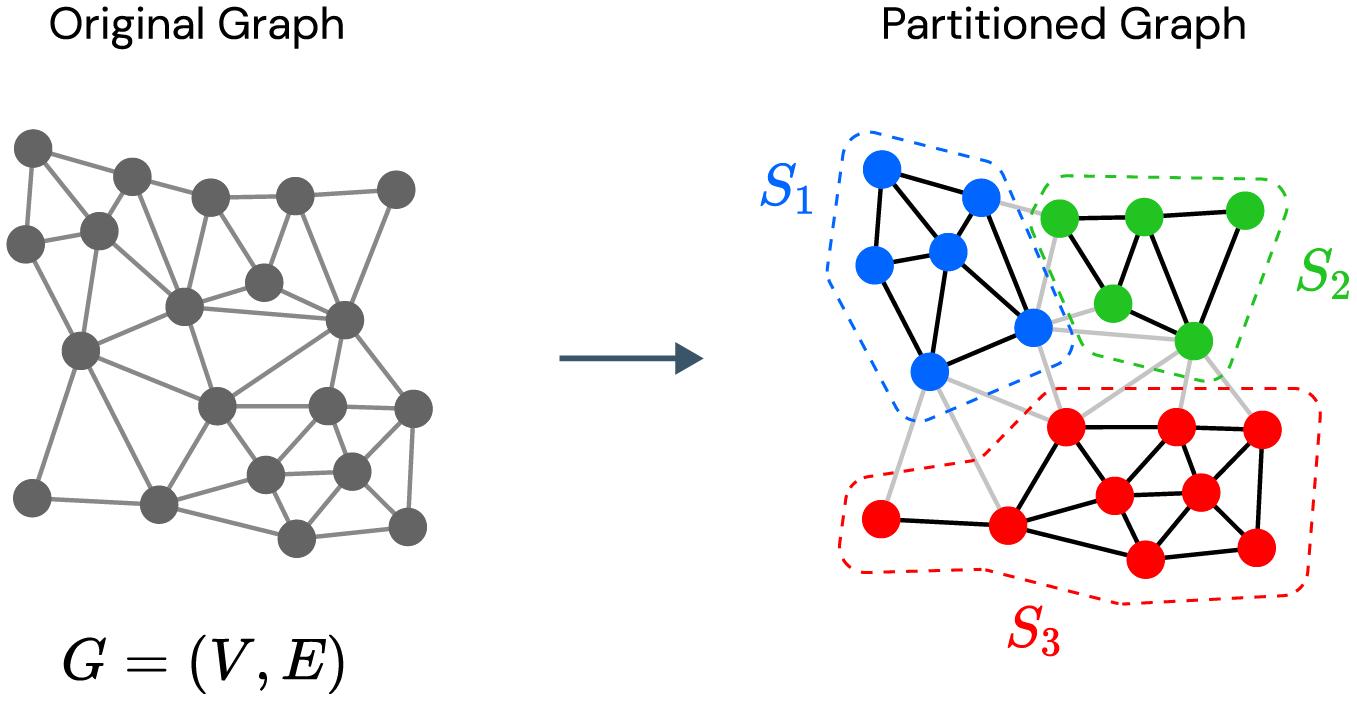}
    \caption{Example of graph partitioning into three sets. Left: original graph. Right: partitioned graph into three subsets $S_1,S_2,S_3$.}
    \label{fig:graphpartitioning}
\end{figure}

Let $N>0$ be the number of graph nodes. Given a graph $G=(V,E)$, where $V=\{v_i\}_{i=1}^N$ and $E=\{e(v_i,v_j):v_i,v_j\in V\}$ are the sets of nodes and the set of edges, respectively. Through this work, we will consider only undirected, unweighted, connected graphs that represent the connectivity of elements of two-dimensional meshes. However, the extension to n-dimensional grids is straightforward, as the connectivity graph is a dimension-less object that can be constructed in principle for meshes of any dimension. Moreover, processing additional information such as directed, weighted edges is trivial for some algorithms, such as the ones based on neural networks, and maybe non-banal or even not possible for others, as we will see in the following.\\
The problem of graph partitioning consists in finding $M$ disjoint sets of nodes $S_1,...,S_M$  such that $\cup_{i=1}^M S_i = V$ and $\cap_{i=1}^M S_i=\emptyset$, see, e.g., Figure \ref{fig:graphpartitioning}.  Sets $S_i,\ i =1,...M$ are assumed to from connected sub-graphs. 
As we restrict our framework to undirected graphs, we have $e(v_i,v_j)\in E \iff e(v_j,v_i)\in E$. We will use the following notation:
\begin{itemize}
    \item $A\in \mathbb{R}^{N\times N}$ is the adjacency matrix, i.e. $A_{i,j} = 1$ if $e(v_i,v_j) \in E$ and 0 otherwise;
    \item $X\in \mathbb{R}^{N\times F}$ is the features matrix of the nodes, which may represent any information such as coordinates, where $F$ is the number of features;
    \item $\mathcal{N}(v_i)=\{v_j\in V: e(v_i,v_j) \in E\}$ denotes the neighborhood of the $i\text{-th}$ node, containing the nodes $v_j$ directly adjacent to $v_i$. In this work, graphs with self-loops are not considered, therefore $v_i \not\in \mathcal{N}(v_i)$.
\end{itemize}
We can define the \textit{cut} of a graph, representing the number of edges connecting the disjoint sets of nodes resulting from the paritioned graph. In the case of two partitions, it can be defined as:
\begin{equation}
    \operatorname{cut}(S_1,S_2) = |\{e(v_i,v_j)\in E : v_i \in S_1, v_j \in S_2\}|, \label{eq:cut}
\end{equation}
and can be easily generalized to the case of $M$ partitions, as
\begin{equation}
    \operatorname{cut}(S_1,...,S_M) = \frac{1}{2}\sum_{k=1}^M \operatorname{cut}(S_k,S_k^C) \qquad \text{with} \quad S_k^C=V\setminus S_k.
\end{equation}
The \textit{normalized cut}  is defined as
\begin{equation}
    \operatorname{Ncut}(S_1,...,S_M) = \sum_{k=1}^M \frac{\operatorname{cut}(S_k,S_k^C)}{\operatorname{vol}(S_k,V)},
\end{equation}
where the \textit{volume} of the partition $S_k$ is defined as
\begin{equation}
     \operatorname{vol}(S_k,V) = |\{e(v_i,v_j) \in E : v_i \in S_k, v_j \in V\}|.
\end{equation}
It represents the total degree of all nodes of the $k\text{-th}$ partition. Depending on the application of interest, different notions of \textit{volume} can be used. The \textit{normalized cut}, contrary to the standard \textit{cut}, allows to take into account partitions where the number of nodes is balanced between the different sets. Let $Y_{ij}$ be the probability for node $i$ of belonging to partition $j$. Considering such a framework allows accounting for the uncertainty of our predictions on a new graph to be partitioned, which some models we will be using can explicitly account for, as we will see in the following. The same framework is reported in \cite{gatti1}. Then, the \textit{expected cut}, given two partitions $S_k$ and its complement $S_k^C$, is defined as
\begin{equation}
    \mathbb{E}[\operatorname{cut}(S_k,S_k^C)] = \sum_{i=1}^N \sum_{v_j \in \mathcal{N}(v_i)} Y_{ik}(1-Y_{jk}) = \sum_{i=1}^N\sum_{j=1}^N Y_{ik}(1-Y_{kj}^T)A_{ij}.
\end{equation}
Let $D$ be the column vector of the degrees of the nodes, i.e. $D_i = \operatorname{vol}(\{v_i\},V)$. Then 
\begin{equation}
    \mathbb{E}[\operatorname{vol}(S_k,V)] = (Y^T D)_k = \Gamma_k.
\end{equation}
The \textit{expected normalized cut} can be defined as
\begin{equation}
\mathbb{E}[\operatorname{Ncut}(S_1,...,S_M)] = \sum (Y\oslash \Gamma)(1-Y)^T\odot A,
\end{equation}
where $\oslash$ and $\odot$ denote the element-wise division and multiplication, respectively, and the summation is over all the entries of the resulting matrix.\\
\newline
A standard algorithm for partitioning large graphs or meshes and computing fill-reducing orderings of sparse matrices is METIS \cite{metis}, based on a multilevel matching procedure, which can also be used as an agglomeration algorithm.\\
Other algorithms for agglomeration that make use of aggregations based on Helmholtz decomposition of the graphs are reported in \cite{hu2021posteriori} and references therein.


\subsection{Machine learning-based graph partitioning}
We propose to employ ML-based bisection models of the form $\mathcal{M}(G,X) = Y$, that take as input a graph $G$ together with a set of features $X$ attached to each node, such as the barycentric coordinates or the area of the mesh elements, and output the vector of probabilities $Y$ of each node belonging to cluster 1 or cluster 2. In order to apply such models for mesh agglomeration, we can use Algorithm~\ref{alg:agglalgo}, that recursively bisect the connectivity graph of the input mesh until the agglomerated elements have the desired size. We recall that the diameter of a domain $\mathcal{D}$ is defined, as usual, as $\textrm{diam}(\mathcal{D}) := \sup \{|x-y|,\ x,y \in \mathcal{D}\}.$ Given a polygonal mesh, i.e. a set of non-overlapping polygonal regions $\mathcal{T}_h=\{P_i\}_{i=1}^{N_P}$, $N_P \geq 1$, that covers a domain $\Omega$, we can define the mesh size $h = \max_{i=1:N_P} \textrm{diam}(P_i).$ Algorithm~\ref{alg:agglalgo} automatically generates a hierarchy of nested grids with different sizes, to be employed e.g. within multigrid solvers. The target meshsize $h^*$ is user-defined and depends on the application. For example, in the context of multigrid solvers a standard choice is to generate grids of size $h_0,\ 2h_0,\ 4h_0$ and so on doubling the size each time, where $h_0$ is the initial meshsize. Adjusting partition $Y$ in line 6 of Algorithm~\ref{alg:agglalgo} is useful in case the provided partition is heuristic and can be therefore sub-optimal. A possible way to adjust the partition is to try minimizing the length of the boundary of the agglomerated mesh elements. This can be done, e.g., by considering each mesh element on the interface of two agglomerated elements, and changing its label accordingly if an improvement in the considered metric is observed.
\begin{algorithm}
\caption{General mesh agglomeration strategy}
\label{alg:agglalgo}
\hspace*{\algorithmicindent} \textbf{Input:} mesh $\mathcal{T}_h$, target mesh size $h^*$, bisection model $\mathcal{M}$. \\
\hspace*{\algorithmicindent} \textbf{Output:} agglomerated mesh $\mathcal{T}_{h^*}$.\\ 

\hspace*{\algorithmicindent} \textbf{Function} \textsc{agglomerate} ($\mathcal{T}_h, \ h^*$)
\begin{algorithmic}[1]
\IF{$\diam(\mathcal{T}_h) \leq h^*$}
\RETURN $\mathcal{T}_h$
\ELSE
\STATE Extract the connectivity graph $G$ and features $X$ from $\mathcal{T}_h$.
\STATE $Y \leftarrow \mathcal{M}(G,X)$
\STATE Adjust partition $Y$.
\STATE Partition $\mathcal{T}_h$ into sub-meshes $\mathcal{T}_h^{(1)},\mathcal{T}_h^{(2)}$ according to $Y$.
\STATE $\mathcal{T}_{h^*}^{(1)} \leftarrow \textsc{agglomerate}(\mathcal{T}_h^{(1)},h^*)$
\STATE $\mathcal{T}_{h^*}^{(2)} \leftarrow \textsc{agglomerate}(\mathcal{T}_h^{(2)},h^*)$
\STATE$\mathcal{T}_{h^*} \leftarrow$ merge   $\mathcal{T}_{h^*}^{(1)}$, $\mathcal{T}_{h^*}^{(2)}$
\ENDIF
\end{algorithmic}
\end{algorithm}
If the model $\mathcal{M}$ does not return a valid partition $Y$, e.g., when the sub-graph within a set is not connected, a suitable fixing procedure is required. This can be done, for example, by considering the largest suitable connected components as new clusters. An efficient algorithm for finding the connected components of a graph (directed or not) is reported in \cite{tarjan1972depth}, which allows to directly consider the blocks of a block-diagonal adjacency matrix.\\
Possible ML-based choices for the bisection model $\mathcal{M}$ are, but not limited to, the k-means clustering algorithm and GNNs. These ML techniques have the advantage to process naturally geometric information about mesh elements, such as their barycentric coordinates, with respect to METIS which is meant to take into account information only about the connectivity graph of the grid. This is important, as it allows to preserve the geometrical quality of the initial grid, enhancing the performance of numerical methods and reducing the number of vertices and edges, therefore reducing also the memory storage and the computational cost associated with the manipulation of the grid \cite{ANTONIETTI2022110900,antonietti2022machine}. For example, starting from an initial grid of squares, it would be desirable that the agglomerated version is also a grid of squares.\\
K-means is a "quantization" algorithm \cite{kmeans,hartigan1979algorithm,likas2003global,bello2012adaptive} that allows partitioning a set of input points in $\mathbb{R}^n$, where $n\geq 1$ is the number of features, into $k$ parts. In order to apply the k-means algorithm for graph bisection within Algorithm \ref{alg:agglalgo}, we employ only two clusters and use as features the barycentric coordinates of the mesh elements. This strategy tends to produce "rounded" (star-shaped) agglomerated elements of similar size. It should be noted that the performance of k-means heavily relies on the initial seeding, which is done as reported in \cite{arthur2007proceedings}. Notice also that no information on the connectivity graph is taken into account directly. However, for a regular mesh with elements of similar size, barycentric coordinates are usually strongly correlated with elements being adjacent.\\
An approach based on GNNs, which we will describe in detail in the following sections, has the following advantages with respect to METIS and k-means:
\begin{itemize}
    \item it can naturally process both information about the connectivity graph of mesh elements and their geometrical features, therefore fully exploiting all of the resources available, while METIS is meant to process only graph-based information and k-means only geometrical information;
    \item it is extremely fast for online inference on new instances, as the main computational burden is faced before offline during a calibration phase known as "training", while METIS and k-means perform iterations online to provide an estimate of the solution; 
    \item it can process additional information in an automatic way, e.g. the  quality of mesh elements according to some metrics, without requiring the user to explicitly know how that information can benefit performance, while for METIS and k-means this has to be explicitly modelled.
\end{itemize}

%% file: sections/GNN.tex
\section{Graph neural networks-based agglomeration strategies}
\begin{figure}
    \centering
    \includegraphics[width=.9\textwidth]{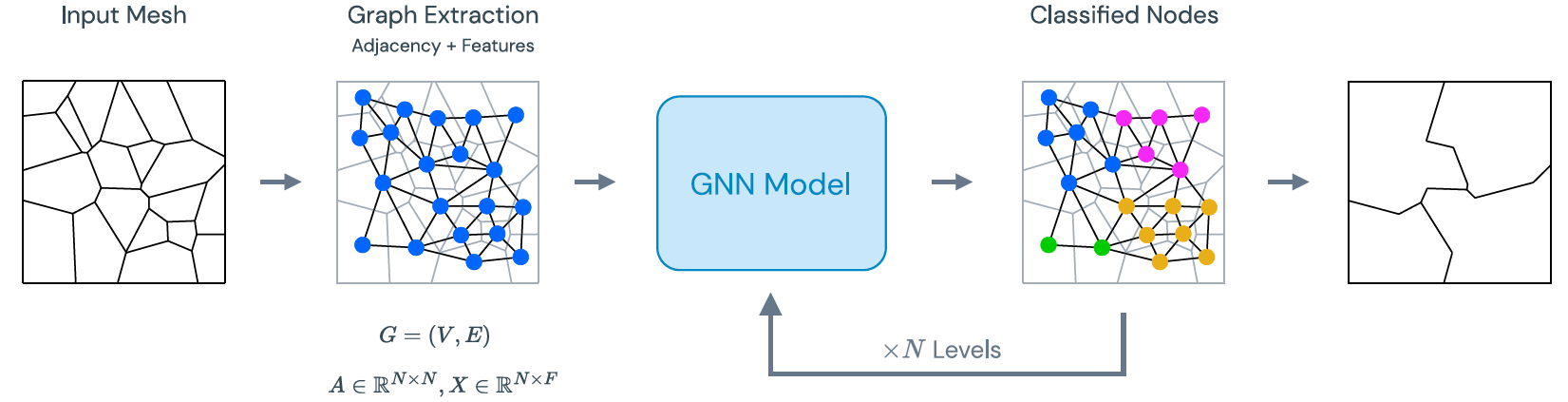}
    \caption{General framework for mesh agglomeration via GNNs.}
    \label{fig:quadtree1}
\end{figure}

In order to perform grid agglomeration by exploiting effectively both the graph representation and the geometrical features of meshes, we employ GNNs as bisection model in Algorithm \ref{alg:agglalgo}. Examples of GNNs-based models directly applied to the original grid are \cite{gatti1,gatti2}, which employ additional processing based on the graph spectrum to perform a preliminary embedding step, in order to extract features that can later be fed to a GNN-based partitioning module. In our case, the graph extraction of geometrical features such the elements area or their barycentric coordinates can be leveraged to perform a classification task, therefore avoiding the need for an additional spectral embedding module. Indeed, spectral bisection algorithms may in general not be well suited for agglomeration purposes as shown in \cite{urschel2016maximal}, where are suitable counter-example is provided. The graph-bisection model performs a classification task, by taking as input the graph $G=(V,E)$ and the features related to its nodes $X\in\mathbb{R}^{N\times F}$, and outputting a probability tensor $Y\in \mathbb{R}^{N\times 2}$, where $Y_{ij}$ represents the probability that the node $v_i \in V$ belongs to the partition $S_j$, with $j=1,2$. This approach can be obviously generalized to an arbitrary number of partitions. However, this would require a specific GNN model for each fixed number of classes, while the 2-classes case can be easily extended to a multi-class case by recursively calling the bisection model on the graphs of each partition. The general framework for mesh agglomeration via GNNs is shown in Figure \ref{fig:quadtree1}.

\subsection{Unsupervised learning for graph partitioning}
In a unsupervised learning framework, we are given an unlabelled dataset, in our case of the form $\{(A^i,X^i)\}_{i=1}^{N_G}$ where $N_G \geq 1$ is the number of graphs in the database, for which we want to find a different representation. For the case of mesh agglomeration, these graphs represent the elements connectivity of different computational grids. We consider then a node classifier $F$, which in our case will be a GNN, parameterized by a set of weights $W$, that takes as input the adjacency matrix $A$ and the nodes features $X$ of a graph and outputs the probability $Y_{ij}$ for node $i$ of belonging to partition $j$. Our goal is to tune $W$ so that $F$ minimizes the following \textit{loss} function
\begin{equation}
    \mathcal{L} = \sum_{i=1}^{N_G} \ell(A^i,Y^i;W),
\end{equation}
where $Y^i = F(A^i,X^i;W)$ and
\begin{equation}
    \ell(A,Y;W) = \sum_{k=1}^M\sum_{i,j=1}^N \frac{Y_{ik}(1-Y_{kj}^T)A_{ij}}{\Gamma_k}
\end{equation}
is the expected normalized cut of the graph, as more extensively addressed in Section \ref{sec: graph partitioning}.

\subsection{Graph neural networks}
\label{sec:gnn_overview}
GNNs are deep learning architectures specifically meant to work with graph-structured data within the framework of Geometrical Deep Learning, that concerns the application of neural networks to non-Euclidean data structures. Mapping, or layers, that can be combined to construct the GNN architecture are the following.

\begin{figure}
    \centering
    \includegraphics[width=1\textwidth]{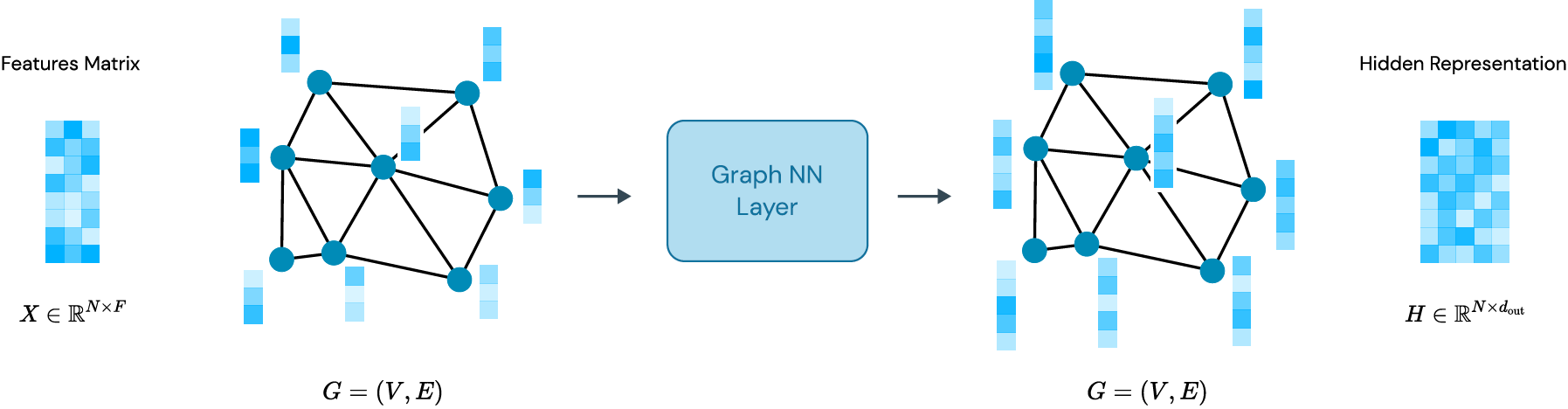}
    \caption{General GNN framework, consisting of an input graph $G$ together with its features matrix $X$, and an output hidden representation $H$ related to the same graph structure.}
    \label{fig:gnnlayer}
\end{figure}

\paragraph{Graph convolutional layers} These layers take as input a graph $G=(V,E)$, consisting of an adjacency matrix $A$ together with features attached to nodes, edges or the global graph, and returns a graph with the same connectivity structure while progressively transforming the information of the features, as shown in Figure \ref{fig:gnnlayer}.
These mappings are permutation-invariant with respect to the order of the nodes. Let $X$ denote the nodes features matrix, or the initial representations, related of the input graph $G$, and let $H^k$ denote the hidden representation of those features after applying the $k\text{-th}$ convolutional layer. At step $k$, the hidden representation $H^k_i$ for the node $v_i$ is computed as:
\begin{align}
    &a^k_i = \Phi^k(\{H_j^{k-1}:v_j\in\mathcal{N}(v_i)\}), \label{eq:aggr} \\
    &H^k_i = \Psi^k(H_i^{k-1},a^k_i), \label{eq:comb} 
\end{align}
where $H^0 = X$ for the initial layer, $\Phi$ is the \textit{aggregation function} that defines how the information coming from the neighborhood $\mathcal{N}(v_i)$ of node $v_i$ is aggregated, and $\Psi$ is the \textit{combination function}  that defines how the aggregated information $a_i^k$ is combined with the one stored in $v_i$. Different definitions of aggregation and combination functions leads to different GNN architectures. In particular, we re-frame equations \eqref{eq:aggr} and \eqref{eq:comb} by considering the mean aggregation function, as follows:
\begin{equation}
\label{eq:SAGEConv}
H_i^{l+1} = \sigma(H_i^l W_1^l + (\operatorname{mean}_{j\in \mathcal{N} (v_i)}H^l_j)W_2^l),
\end{equation}
where $\sigma(\cdot)$ is a non-linear activation function, such as the REctified Linear Unit (ReLU) or the hyperbolic tangent ($tanh$), and $W_1^k,W_2^k\in \mathbb{R}^{F_k\times F_{k+1}}$ are weight matrices associated to the $l\text{-th}$ layer, representing a trainable linear transformation, where $F_k$ and $F_{k+1}$ are the features dimensions for the current and next layers, respectively. We refer to \eqref{eq:SAGEConv} as SAmpling-and-aggreGatE Convolutional (\textsc{SAGEConv}) layer \cite{1706.02216} followed by activation function $\sigma$. Such layers also account for the possibility of sub-sampling the neighbourhood of a node during the aggregation process, leveraging the information coming from features to better generalize to unseen nodes, while keeping the computational cost under control.

    \paragraph{Input normalization layer} We consider a mapping of the form $\textsc{INorm}: \mathbb{R}^{N \times F} \rightarrow \mathbb{R}^{N \times F}$, where $N$ is the number of nodes and $F$ is the number of features, that normalizes the input feature matrix $X = [x_1 | \dots | x_F]$. The normalization is performed differently for each type of feature, i.e. column-wise: for strictly positive features, such as the mesh elements areas, we simply rescale to $[0,1]$
    \begin{equation*}
        \tilde{x}_i = \frac{x_i}{\max(x_i)}, \quad i = 1,..., F,
    \end{equation*}
    while for other features, such as barycentric coordinates, we first center them, by subtracting the mean, and then rescale to $[-1,1]$
    $$
         \tilde{x}_i = \frac{y_i}{\max(|y_i|)}, \quad y_i = x_i-\operatorname{mean}(x_i), \quad i = 1,..., F.
    $$
    
    \paragraph{Dense layers} Dense layers are generic linear maps of the form $\textsc{Linear}:  \mathbb{R}^{m} \rightarrow \mathbb{R}^\ell,$ $m,\ell \geq 1$ defined by parameters to be tuned. They are used to separate graph features extracted in the previous layers.
    
    \paragraph{Softmax} We define the function $\textsc{Softmax}: \mathbb{R}^\ell \rightarrow (0,1)^\ell,$ where $\ell\geq 2$ is the number output classes, $$[\textsc{Softmax}(x)]_i = \frac{e^{x_i}}{\sum_{j=1}^\ell e^{x_j}}.$$ They are used to assign a probability to each class.\\[0.5em]

\subsection{Graph neural network training}
\begin{figure}
    \centering
    \includegraphics[width=1\textwidth]{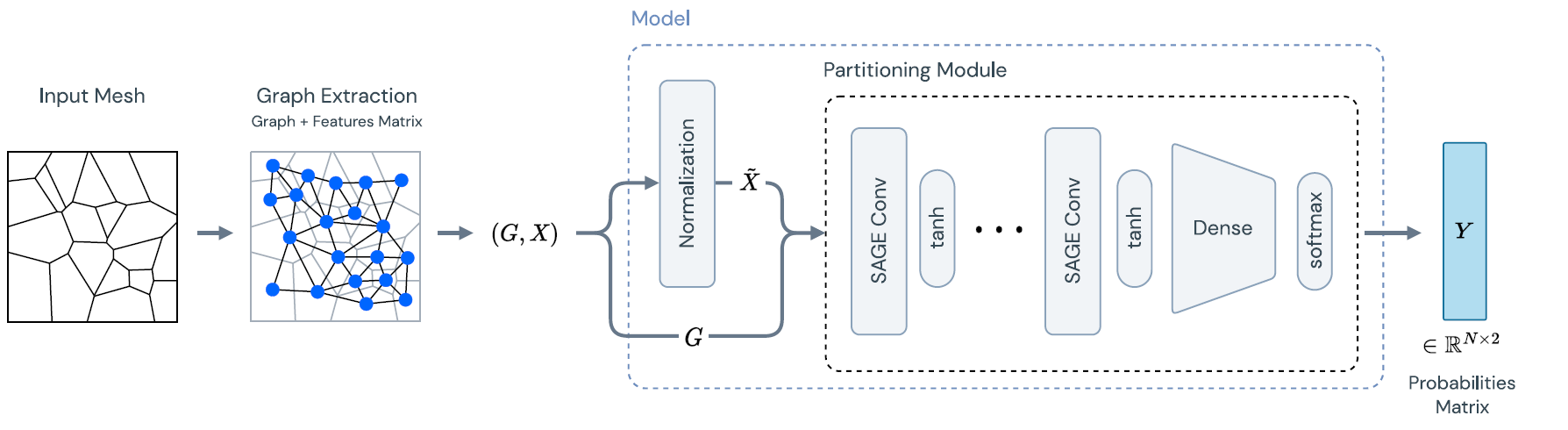}
    \caption[GNN model structure.]{GNN model structure, consisting of a normalization module, responsible of rescaling the features extracted from the input mesh, and of a partitioning module made of a stack of graph convolutions and a dense classifier. The output is a matrix of probabilities $Y\in\mathbb{R}^{N\times 2}$, containing the probabilities of belonging either to one of the two classes for each node.}
    \label{fig:quadtree2}
\end{figure}

The input features matrix is $X \in \mathbb{R}^{N\times 3}$, where the first column contains the area of each mesh element followed by the two coordinates of its barycenter. The GNN architecture we employed first applies a \textsc{INorm} layer, then four \textsc{SAGEConv} layer with features dimensions 64, then three \textsc{Linear} layer with progressively decreasing output dimensions (32, 8 and 2) and finally a \textsc{Softmax} layer. Each \textsc{SAGEConv} layer is followed by a $tanh$ activation function to keep the features inside the interval $[-1,1]$, so that the geometrical information regarding the re-scaled coordinates will be kept in the same domain, as information flows through the layers. To further simplify the classification process, the \textsc{INorm} layer also rotates of 90 degrees the barycentric coordinates if the input mesh is more stretched along the y-axis rather than the x-axis. The resulting model consists of approximately 28k parameters, where roughly 25k resulting from the \textsc{SAGEConv} layers and the remaining from the \textsc{Linear} layers. A scheme of the described GNN model is shown in Figure \ref{fig:quadtree2}.\\
\newline
The training of the models has been performed by employing an unsupervised approach by minimizing the expected normalized cut, as described in Section~\ref{sec:graphpart}. To perform the training a train and a validation datasets were generated. The training dataset consists of meshes of the following type: grids of regular squares, grids of regular triangles, grids of triangles with random perturbation of the vertices, grids of Voronoi with random location of the seeds. The database is composed of 800 meshes, with 200 meshes per type, while the validation dataset consists of 200 meshes, 50 per type. The cardinality of the datasets has been chosen to keep a 80-20 split ratio between train and validation respectively. Training has been performed using the Adam optimizer \cite{kingma2015adam} with a learning rate 1e-5, $L^2$ regularization coefficient 1e-5 and mini-batch size 4. Training was performed for 300 epochs in approximately 35 minutes on Google Colab cloud platform, using a 2.20GHz Intel Xeon processor, with 12GB of RAM memory and NVIDIA Tesla T4 GPU with 16GB of GDDR6 memory.

%% file: sections/validation_grids.tex
\section{Validation on a set of polyhedral grids}\label{sec:results}
\begin{figure}
    \centering
    \includegraphics[width=\textwidth]{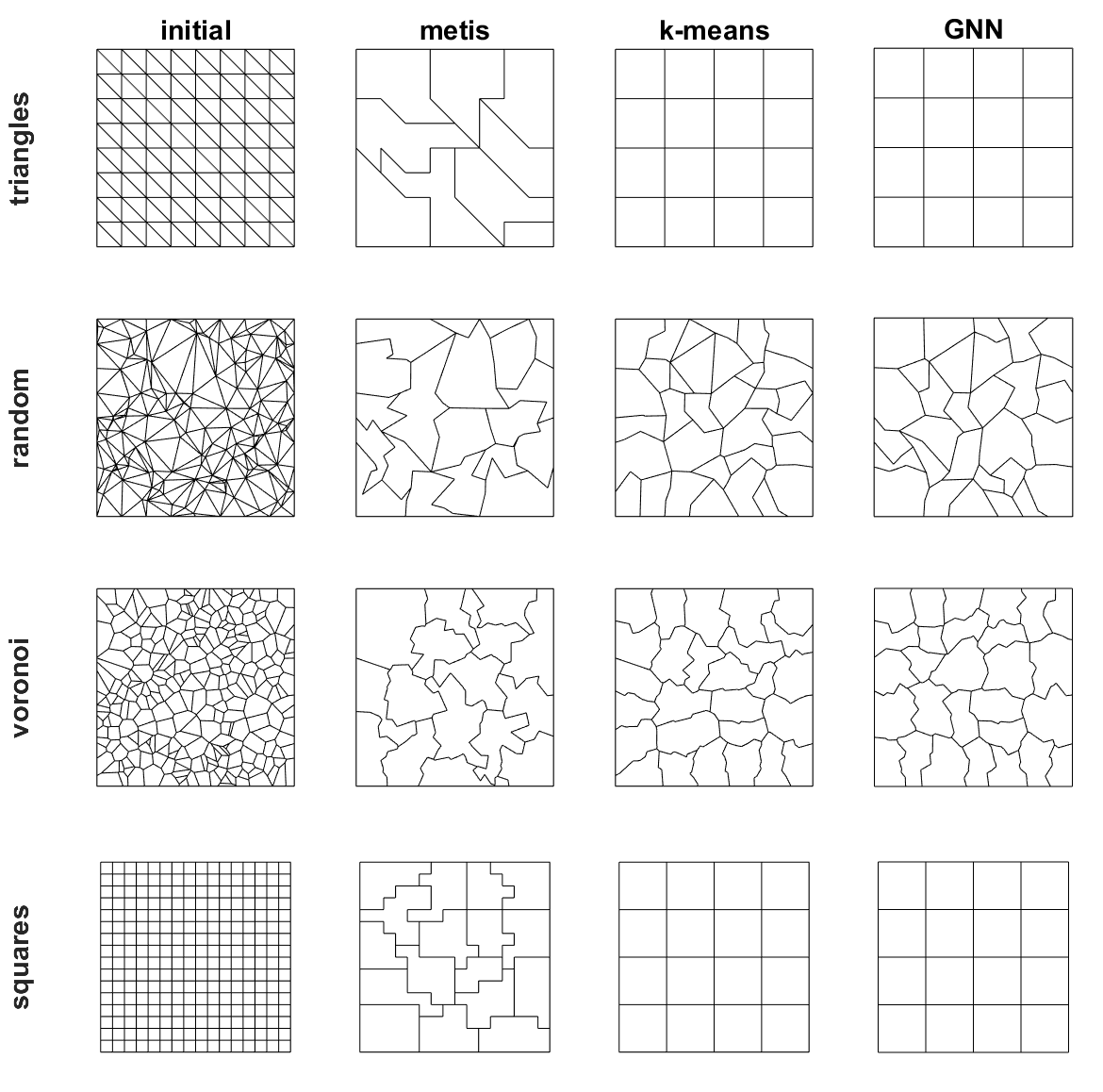}
    \caption[Initial grids agglomerated using different agglomeration strategies.]{Initial grids (first column) and corresponding agglomerated versions (second to fourth column) obtained based on employing different strategies. Each row corresponds to the same initial grid (from top to bottom: regular triangles, random triangles, Voronoi, squares) while each column corresponds to the same agglomeration strategy (from left to right: initial grids, METIS, k-means, GNN).}
    \label{fig:uniformaggl}
\end{figure}

\begin{figure}
    \centering
    \includegraphics[width=.9\textwidth]{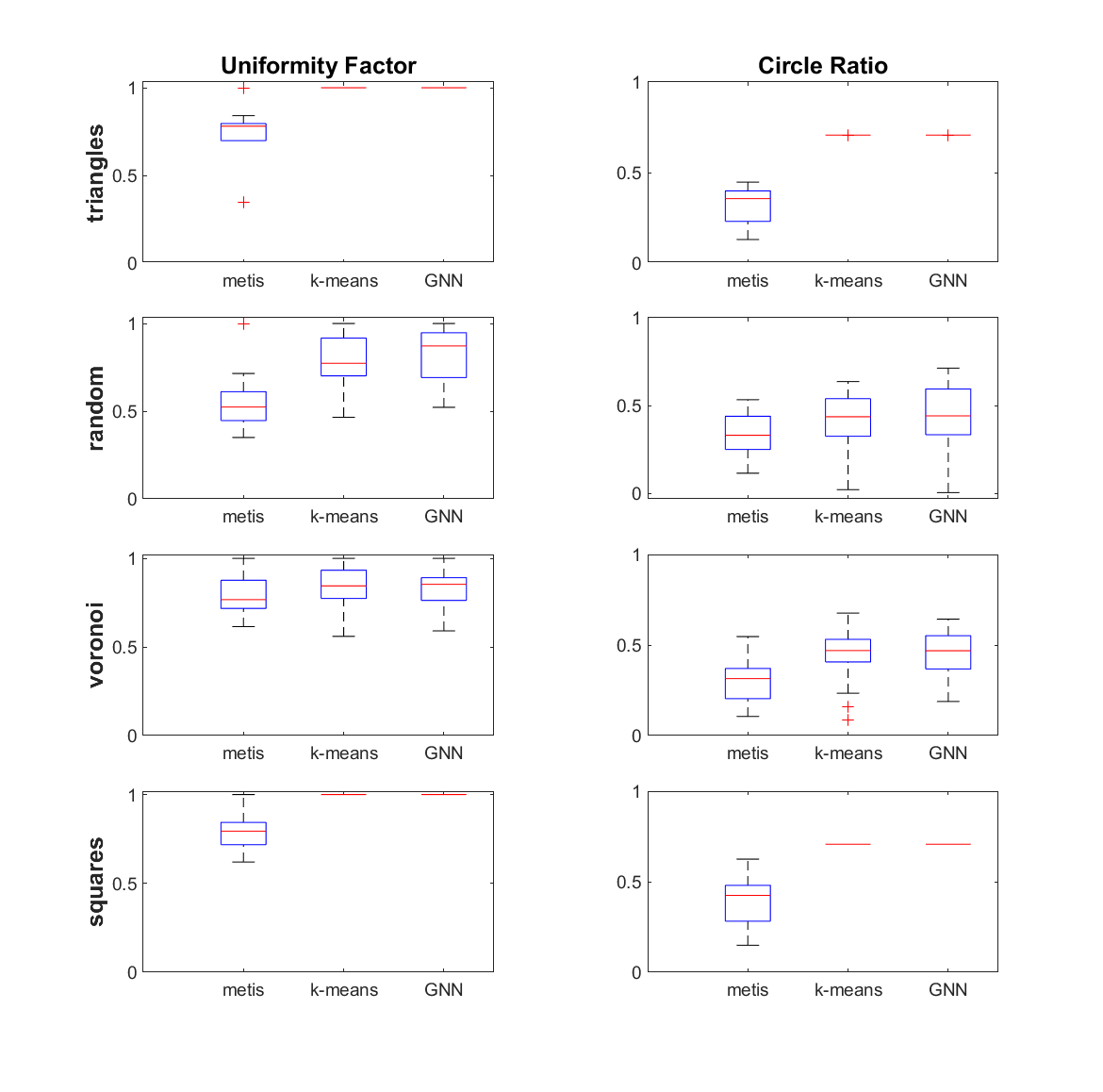}
    \caption[Box plots of the computed quality metrics (UF and CR) for the agglomerated grids.]{Box plots of the computed quality metrics (UF and CR) for the agglomerated grids reported in Figure~\ref{fig:uniformaggl} (second to fourth column), obtained based on employing different agglomeration strategies (METIS, k-means, GNN). Some box plots collapse into a single line because the corresponding metric has the same value for all mesh elements.}
    \label{fig:qualityplot}
\end{figure}

In this section we compare the performance of the proposed algorithms. To evaluate the quality of the refined grids, we employ the following quality metrics introduced in \cite{attene2019benchmark}:
\begin{itemize}
    \item \textit{Uniformity Factor} (UF): ratio between the diameter of an element $P$ and the mesh size
    $$\textrm{UF}(P) = \frac{\textrm{diam}(P)}{h}.$$
    This metric takes values in $[0,1]$. The higher its value is the more mesh elements have comparable sizes.
    \item \textit{Circle Ratio} (CR): ratio between the radius of the inscribed circle and the radius of the circumscribed circle of a polyhedron $P$
    $$\textrm{CR}(P) = \frac{\max_{\{B(r) \subset P\}}  r }{\min_{\{P \subset B(r)\}} r },$$
	where $B(r)$ is a circle of radius $r$. For the practical purpose of measuring the roundness of an element the radius of the circumscribed circle has been approximated with $\textrm{diam}(P) / 2$. This metric takes values in $[0,1]$. The higher its value is the more rounded mesh elements are.
\end{itemize}
We consider four different grids of domain $(0, 1)^2$: a grid of regular triangles, a grid of triangles with random location of the vertices, a Voronoi grid with random location of the seeds, and a grid of regular squares. In Figure \ref{fig:uniformaggl} these grids have been agglomerated using METIS, k-means and GNNs strategies. For METIS the target number of mesh elements is $N_0/16$, where $N_0$ is the initial number of elements, while for k-means and GNN the target mesh size in Algorithm \ref{alg:agglalgo} is $4h_0$, where $h_0$ is the initial mesh size. 
As we can see, the k-means and the GNN algorithms are capable to recover a regular grid of squares when starting from regular meshes (triangles and squares), while this is not the case for METIS. In Figure~\ref{fig:qualityplot} we show the box plots of the computed quality metrics for the selected grids.
In general, quality metrics are lower for the METIS algorithm, while k-means and GNNs have comparable performance. This is further confirmed by Tables \ref{table:UF} and \ref{table:CR}, where we report the average values of UF and CR. In particular, the performance difference is much more evident for regular grids. 
\begin{table}
\centering
\begin{tabular}{|c|ccc|}
\hline
UF        & metis  & k-means & GNN    \\ \hline
triangles & 0.7648 & 1.0000  & 1.0000 \\
random    & 0.6115 & 0.7003  & 0.7418 \\
Voronoi   & 0.7605 & 0.8105  & 0.8225 \\
squares   & 0.7725 & 1.0000  & 1.0000 \\ \hline
\end{tabular}
\caption[Average values of the Uniformity Factor for the agglomerated grids.]{Average values of the UF for the agglomerated grids reported in Figure~\ref{fig:uniformaggl}, obtained based on employing different agglomeration strategies (METIS, k-means, GNN).}
\label{table:UF}
\end{table}
\begin{table}
\centering
\begin{tabular}{|c|ccc|}
\hline
CR        & metis  & k-means & GNN    \\ \hline
triangles & 0.3447 & 0.7071  & 0.7071 \\
random    & 0.3399 & 0.4190  & 0.3928 \\
Voronoi   & 0.3056 & 0.4509  & 0.4845 \\
squares   & 0.3733 & 0.7071  & 0.7071 \\ \hline
\end{tabular}
\caption[Average values of the Circle Ratio for the agglomerated grids.]{Average values of the CR for the agglomerated grids reported in Figure~\ref{fig:uniformaggl}, obtained based on employing different agglomeration strategies (METIS, k-means, GNN).}
\label{table:CR}
\end{table}
In Tables \ref{table:UF rel} and \ref{table:CR rel} we also report the relative performance with respect to METIS, i.e., the ratio between the average UF and CR metrics of the ML-strategies (k-means and GNN) and METIS. This further highlights the gain is using ML-based strategies.
\begin{table}
\centering
\begin{tabular}{|c|cc|}
\hline
UF relative & k-means/metis & GNN/metis    \\ \hline
triangles   & 1.3076 & 1.3076 \\
random      & 1.1452 & 1.2130 \\
Voronoi     & 1.0658 & 1.0815 \\
squares     & 1.2945 & 1.2945 \\ \hline
\end{tabular}
\caption[Relative Uniformity Factor performance with respect to METIS for the agglomerated grids.]{Relative UF: performance improvement with respect to METIS, i.e., ratio between the average UF of the ML-strategies (k-means and GNN) and METIS, for the agglomerated grids reported in Figure~\ref{fig:uniformaggl}.}
\label{table:UF rel}
\end{table}
\begin{table}
\centering
\begin{tabular}{|c|cc|}
\hline
CR relative & k-means/metis & GNN/metis    \\ \hline
triangles   & 2.0512  & 2.0512 \\
random      & 1.2327  & 1.1554 \\
Voronoi     & 1.4755  & 1.5853 \\
squares     & 1.8940  & 1.8940 \\ \hline
\end{tabular}
\caption[Relative Circle Ratio performance with respect to METIS for the agglomerated grids.]{Relative CR: performance improvement with respect to METIS, i.e., ratio between the average CR of the ML-strategies (k-means and GNN) and METIS, for the agglomerated grids reported in Figure~\ref{fig:uniformaggl}.}
\label{table:CR rel}
\end{table}
In general, ML-based strategies (either the k-means and the GNN algorithms), seem to preserve the initial geometry and quality of the grids, because the geometric information attached to the nodes is taken into account, making them suitable for adaptive mesh coarsening. This is not the case for the METIS algorithm, because it processes only the information coming from the graph topology of the mesh.

\subsection{Application to a computational mesh stemming from a human brain MRI-scan}
\begin{figure}
    \centering
    \includegraphics[width = \linewidth]{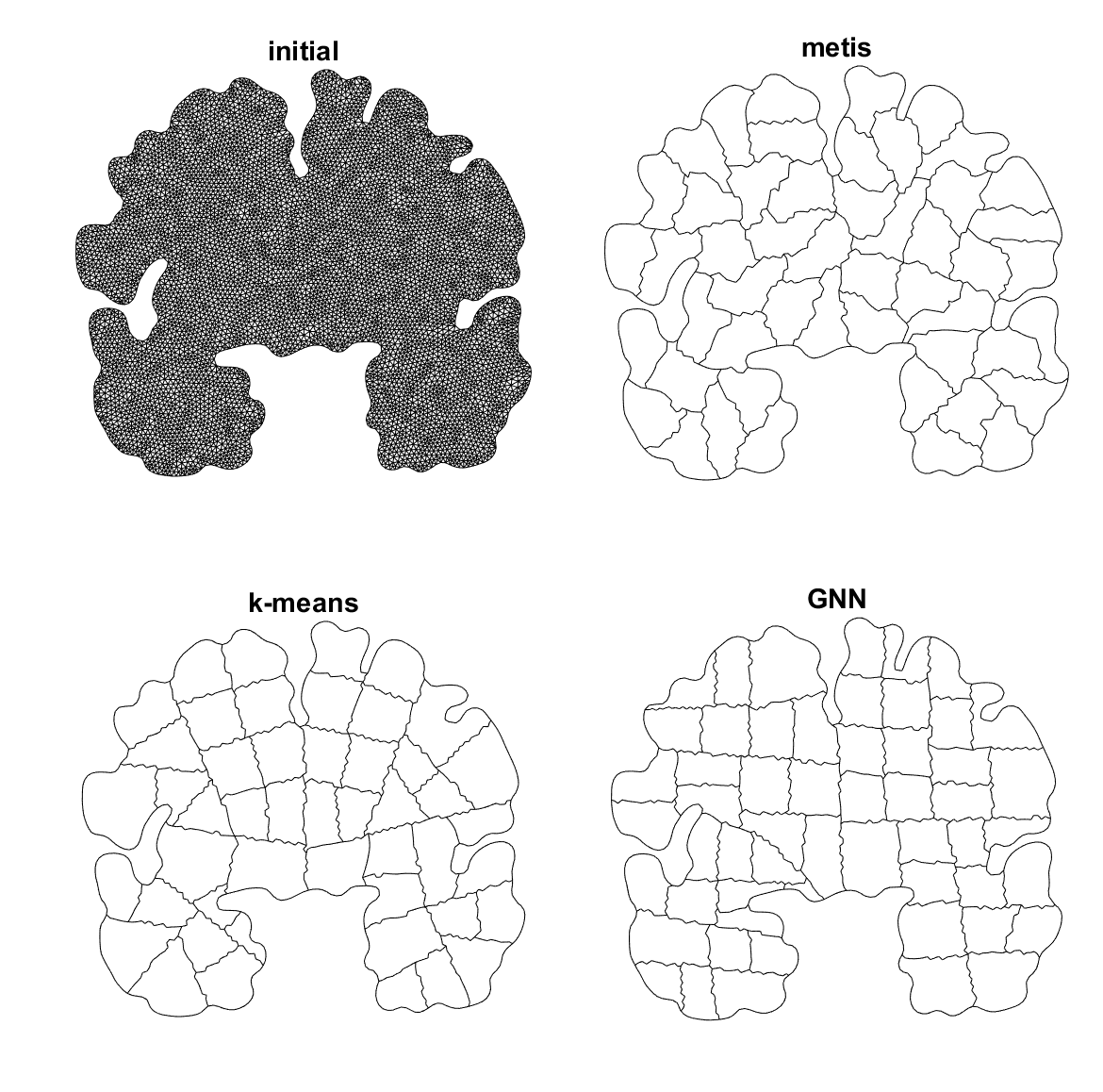}
    \caption[Agglomerated meshes using different strategies, starting from an initial grid of a human brain MRI-scan.]{Agglomerated meshes using different strategies (METIS, k-means, GNN), starting from an initial grid of a human brain MRI-scan, consisting of 14372 triangular elements as shown in the top-left figure.}
    \label{fig:brainaggl}
\end{figure}

\begin{figure}
    \centering
    \includegraphics[width=1\textwidth]{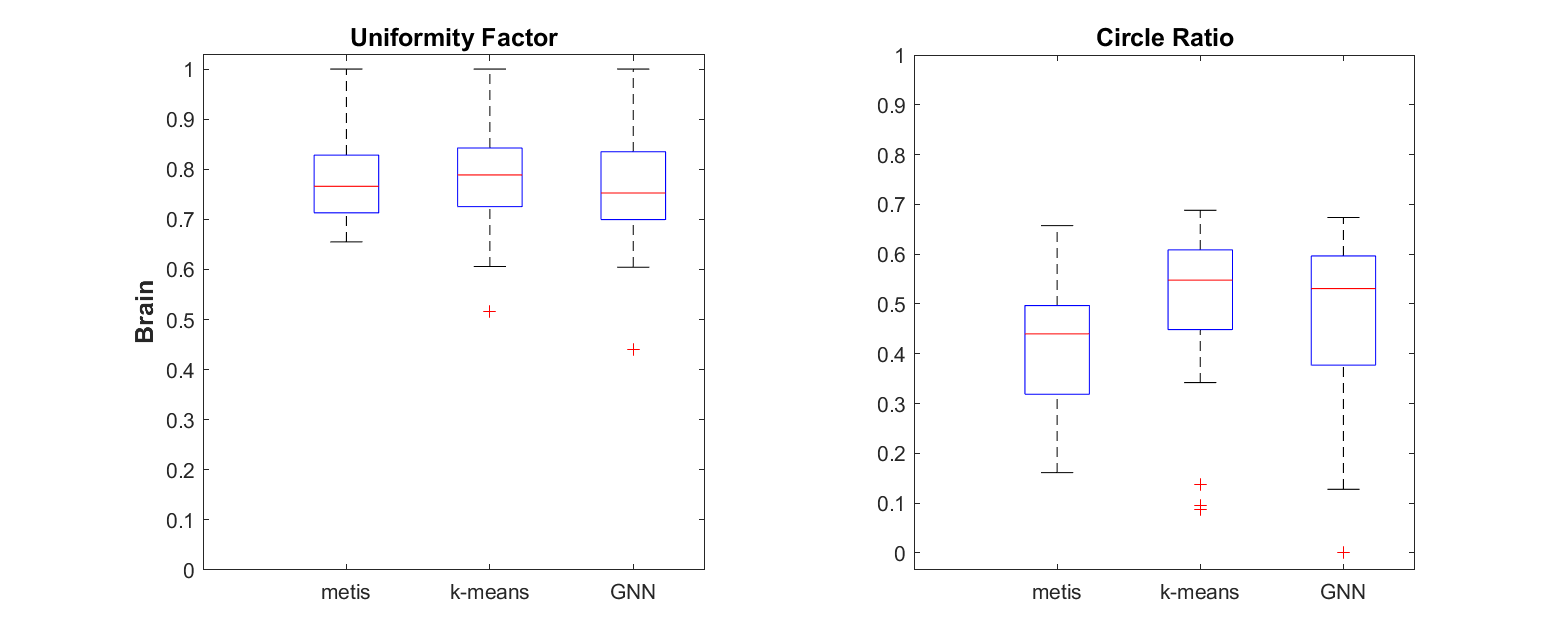}
    \caption[Box plots of the computed quality metrics (Uniformity Factor and Circle Ratio) for the agglomerated MRI brain scan meshes.]{Box plots of the computed quality metrics (UF and CR) for the agglomerated MRI brain scan meshes reported in Figure~\ref{fig:brainaggl}, obtained based on employing different agglomeration strategies (METIS, k-means, GNN).}
    \label{fig:qualitybrain}
\end{figure}

In order to further test the generalization capabilities of our models, we apply them on a much more complex domain with respect to the meshes considered so far. In particular, we consider the mesh of a section of a human brain coming from an MRI-scan, consisting of 14372 triangular elements. The domain is highly non-convex and presents many constrictions and narrowed sections. We agglomerated such a grid using
METIS, k-means and GNN algorithms: the result is reported in Figure \ref{fig:brainaggl}. For k-means and GNN the target mesh size in Algorithm \ref{alg:agglalgo} is $0.2D$, where $D$ is the maximum distance between any two vertices of the initial mesh. For METIS the target number of mesh elements is 50, which corresponds approximately to the same mesh size.
In Figure \ref{fig:qualitybrain} we show also the corresponding box plots of the quality metrics and in Table~\ref{table:metricsbrain} we report the average values.
\begin{table}
\centering
\begin{tabular}{@{}|c|ccc|@{}}
\hline
metric      & metis  & k-means & GNN    \\ \hline
UF          & 0.7808 & 0.7875  & 0.7672 \\
CR          & 0.4134 & 0.5146  & 0.4841 \\
UF relative & 1      & 1.0085  & 0.9826 \\
CR relative & 1      & 1.2449  & 1.1712 \\ \hline
\end{tabular}
\caption[Average and relative Uniformity Factor and the Circle Ratio for agglomerated MRI brain scan mesh.]{Average and relative values (with respect to METIS) of the UF and the CR for different agglomeration methods (METIS, k-means, GNN), applied to the section of the MRI brain scan, for the agglomerated MRI brain scan meshes reported in Figure~\ref{fig:brainaggl}.}
\label{table:metricsbrain}
\end{table}
As we can see, performance are comparable in terms of UF, while the ML-based strategies achieve on average a higher "roundness", measured in terms of CR. This indicates good generalization capabilities of the GNN algorithm, as the considered mesh was very different from the ones included in the training set, both in terms of shape, dimensions and number of elements. The k-means reasonably performs slightly better than the GNN, because it does not require prior information coming from a database at the cost of having a higher online computational cost, as we will see in the next section. 


\subsection{Runtime performance}
\label{Perform}
\begin{figure}
    \centering
    \includegraphics[width = \linewidth]{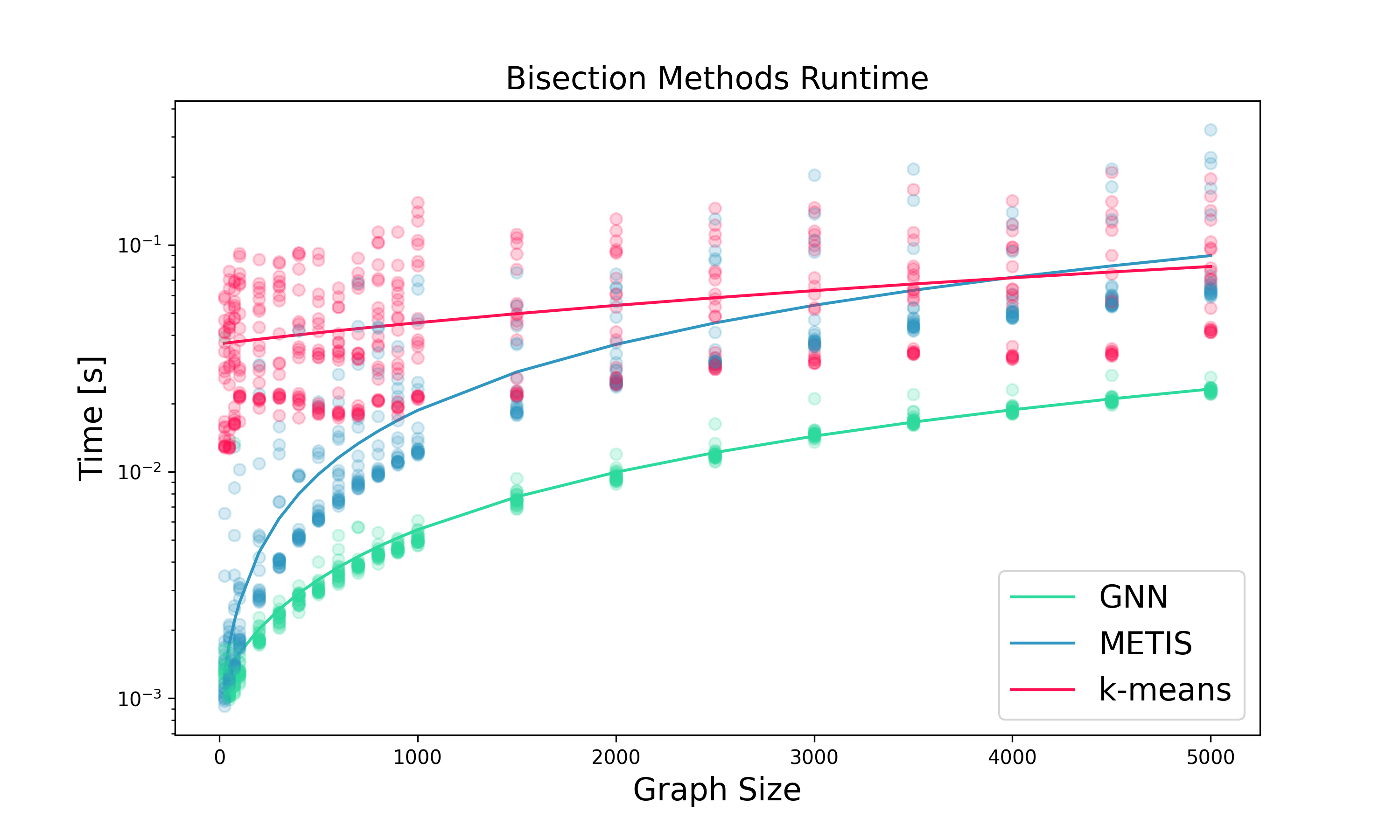}
    \caption[Runtime performance for different graph bisection models.]{Runtime performance for different graph bisection models (METIS, k-means, GNN) as a function of the number of nodes in the connectivity graph of Voronoi meshes. The y-axis is in logarithmic scale.}
    \label{fig:runtimeplot}
\end{figure}

It is well known already that the graph bisection problem is NP-hard \cite{gary1979computers,andreev2004balanced}, meaning that it is solved approximately with heuristic techniques. In general, it is hard to establish how far these approximations are from the optimal solution, and performance may depend heavily on the specific instance of the problem which is being considered.
The GNN, k-means and METIS algorithms all have a computational complexity that scales linearly with respect to the number of mesh elements, and all of them can be parallelized \cite{karypis1997parmetis,kantabutra2000parallel,ma2018towards}.
However, the main advantage of using an approach based on neural networks is the low computational cost for online inference of a new problem instance, as the main cost is offline during training. Moreover, the size of the network mainly scales with the dimension of the data manifold \cite{petersen2020neural}, and therefore also the associated complexity. On the contrary, the computational cost associated with METIS and k-means is online, as these methods need to perform iterations, respectively the coarsening-uncorseaning steps for METIS and the Lloyd's iterations for k-means, in order to improve the current estimate.\\
In order to measure the online performance of the three methods, we applied them these graph partitioning algorithms on 21 random Voronoi meshes with an increasing number of elements from 25 to 5000. Since the runtimes are not deterministic, we sampled 20 runtimes for each mesh. Results are shown in Figure \ref{fig:runtimeplot}.
As expected, the GNN algorithm outperforms the METIS and the k-means ones. For a mesh with 5000 elements we have that GNN is 4.4836 times faster than Metis and 3.5670 times faster than k-means. Being able to compute the solution in a fast manner, even with a slight loss of accuracy, can be also particularly beneficial when using multigrid schemes as preconditioners.

%% file: sections/multigrid.tex
\section{Applications to agglomeration-based multigrid methods}

In this section we test the effectiveness of the proposed agglomeration strategies, to be used in combination with polygonal finite element discretizations within a MultiGrid (MG) framework \cite{bassi2012flexibility,bassi2012agglomeration,antonietti2015multigrid,antonietti2017multigrid,xu2017algebraic,chan1998multilevel,multigridV}. We consider the following model problem: find $u \in V = \text{H}^2(\Omega) \cap \text{H}^1_0(\Omega)$ such that
\begin{equation}
 \label{eq:Poisson problem MG}
\int_\Omega \nabla u \cdot \nabla v = \int_\Omega f v \quad \forall v \in V,
\end{equation}
with $\Omega$ and $f \in L^2(\Omega)$ the forcing term, selected in such a way that the exact solution is given by $u(x,y) = \sin(\pi x)\sin(\pi y).$
We consider the V-cycle MG algorithm with additive Schwarz smoothing within a PolyDG discretization framework, as described in \cite{multigridV}. We employ the initial grids shown in the first column of Figure \ref{fig:uniformaggl} (triangles, random, Voronoi, squares), agglomerated using the proposed strategies (METIS, k-means, GNNs). In particular, for each initial mesh we construct three agglomerated grids of increasing size: for METIS the target numbers of mesh elements are $N_0/4$, $N_0/16$, $N_0/64$, where $N_0$ is the initial number of elements, while for k-means and GNN the target mesh sizes in Algorithm \ref{alg:agglalgo} are $2h_0$, $4h_0$, $8h_0$, where $h_0$ is the initial mesh size. These correspond to the different levels of the V-cycle algorithm, where level 1 corresponds to the initial grid (level 3 grids are the ones reported in Figure \ref{fig:uniformaggl}). As a performance metric, we consider the iteration counts of the MG algorithm to reduce the (relative) residual below $10^{-6}$ in solving the algebraic formulation of problem \eqref{eq:Poisson problem MG}. We vary the following parameters: number of levels employed $\ell$, polynomial degree $p$, number of smoothing steps $m$. In Table \ref{table:MG L m3} we report the iteration counts when varying the number of levels $\ell = 2,3,4$ with $m = 3$ and $p = 1$. As a comparison, the iteration counts of the Conjugate Gradient (CG) method are also reported.


\begin{table}
\centering
\begin{tabular}{|c|c|ccc|c|}
\hline
\multirow{2}{*}{grids} & \multirow{2}{*}{$\ell$} & \multicolumn{3}{l|}{Agglomeration-based MG} & \multirow{2}{*}{CG} \\ \cline{3-5}
                       &                             & metis         & k-means        & GNN        &                     \\ \hline
          & 2      & 21    & 9       & 9   &     \\
triangles & 3      & 21    & 9       & 9   & 114 \\
          & 4      & 21    & 9       & 9   &     \\ \hline
          & 2      & 46    & 41      & 32  &     \\
random    & 3      & 46    & 41      & 32  & 655 \\
          & 4      & 46    & 41      & 32  &     \\ \hline
          & 2      & 19    & 16      & 15  &     \\
Voronoi   & 3      & 19    & 16      & 15  & 348 \\
          & 4      & 19    & 16      & 15  &     \\ \hline
          & 2      & 20    & 17      & 17  &     \\
squares   & 3      & 20    & 17      & 17  & 109 \\
          & 4      & 20    & 17      & 17  &     \\ \hline
\end{tabular}
\caption[Multigrid iteration counts with different levels and three smoothing steps.]{iteration counts of the MG algorithm to reduce the (relative) residual below $10^{-6}$ employing different initial grids (triangles, random, Voronoi, squares) agglomerated with different strategies (METIS, k-means, GNN) with $\ell = 2,3,4$, $p = 1$, $m = 3$. As a comparison, the iteration counts of the Conjugate Gradient (CG) method are reported in the last column.}
\label{table:MG L m3}
\end{table}
We can observe the following:
\begin{itemize}
    \item The iteration counts of the MG methods are significantly lower than the ones of the classical CG method. This is expected, as the CG method only employs the finest grid available without the need to resort to agglomeration techniques. This confirms that the considered MG implementations are indeed effective.
    \item The iteration counts of the k-means and GNN algorithms are lower than the ones of METIS, meaning the higher grid quality provided by the ML-based agglomeration strategies can help accelerating the convergence of the numerical method, with GNN having the best performance.
    \item The number of iterations required when varying the number of levels $\ell$ is constant, meaning it is independent of the granularity of the underlying grid and therefore scalable in terms of the mesh size.
\end{itemize}
In Table \ref{table:MG L m1} we also consider the case for $m = 1$.
\begin{table}
\centering
\begin{tabular}{|c|c|ccc|c|}
\hline
\multirow{2}{*}{grids} & \multirow{2}{*}{$\ell$} & \multicolumn{3}{l|}{Agglomeration-based MG} & \multirow{2}{*}{CG} \\ \cline{3-5}
                       &                             & metis         & k-means        & GNN        &                     \\ \hline
\multirow{2}{*}{triangles} & 3      & 230   & 67      & 67  & \multirow{2}{*}{114} \\
          & 4      & 227   & 66      & 66  &     \\ \hline
\multirow{2}{*}{random }   & 3      & 706   & 551     & 405 & \multirow{2}{*}{655} \\
          & 4      & 708   & 577     & 407 &     \\ \hline
\multirow{2}{*}{Voronoi}   & 3      & 154   & 143     & 129 & \multirow{2}{*}{348} \\
          & 4      & 154   & 143     & 131 &     \\ \hline
\multirow{2}{*}{squares}   & 3      & 178   & 162     & 162 & \multirow{2}{*}{109} \\
          & 4      & 180   & 162     & 162 &     \\ \hline
\end{tabular}
\caption[Multigrid iteration counts with different levels and one smoothing step.]{iteration counts of the MG algorithm to reduce the (relative) residual below $10^{-6}$ employing different initial grids (triangles, random, Voronoi, squares) agglomerated with different strategies (METIS, k-means, GNN) with $\ell = 3,4$, $p = 1$, $m = 1$. As a comparison, the iteration counts of the Conjugate Gradient (CG) method are reported in the last column.}
\label{table:MG L m1}
\end{table}
In order to guarantee a proper convergence of the MG method, the number of smoothing steps should be sufficiently large, as reported in \cite{multigridV}. Despite using such a low value, the MG method still converges. As expected, the iteration counts increase but the ones of the ML strategies are still lower with respect to the ones of the CG, while this is not the case for the METIS algorithm.\\ 
In Table \ref{table:MG m} we report the iteration counts when varying the number of smoothing steps $m = 3,5$ with $\ell = 3$ and $p = 1$.
\begin{table}
\centering
\begin{tabular}{|c|c|ccc|c|}
\hline
\multirow{2}{*}{grids} & \multirow{2}{*}{\textit{m}} & \multicolumn{3}{l|}{Agglomeration-based MG} & \multirow{2}{*}{CG} \\ \cline{3-5}
                       &                             & metis         & k-means        & GNN        &                     \\ \hline
\multirow{2}{*}{triangles} & 3          & 21    & 9       & 9   & \multirow{2}{*}{114} \\
          & 5          & 10    & 5       & 5   &     \\ \hline
\multirow{2}{*}{random}   & 3          & 46    & 41      & 32  & \multirow{2}{*}{655} \\
          & 5          & 20    & 18      & 15  &     \\ \hline
\multirow{2}{*}{Voronoi}   & 3          & 19    & 16      & 15  & \multirow{2}{*}{348} \\
          & 5          & 10    & 8       & 8   &     \\ \hline
\multirow{2}{*}{squares}   & 3          & 20    & 17      & 17  & \multirow{2}{*}{109} \\
          & 5          & 9     & 8       & 8   &     \\ \hline
\end{tabular}
\caption[Multigrid iteration counts with different smoothing steps.]{iteration counts of the MG algorithm to reduce the (relative) residual below $10^{-6}$ employing different initial grids (triangles, random, Voronoi, squares) agglomerated with different strategies (METIS, k-means, GNN) with $\ell = 3$, $p = 1$, $m = 3,5$. As a comparison, the iteration counts of the Conjugate Gradient (CG) method are reported in the last column.}
\label{table:MG m}
\end{table}
As we can see, when the number of steps increases the iteration counts decrease, as well as the performance difference between the different methods. In Table \ref{table:MG p} we report the iteration counts when varying the polynomial degree $p = 1,2,3$ with $\ell = 3$ and $m = 3$.
\begin{table}
\centering
\begin{tabular}{|c|c|ccc|c|}
\hline
\multirow{2}{*}{grids} & \multirow{2}{*}{\textit{p}} & \multicolumn{3}{l|}{Agglomeration-based MG} & \multirow{2}{*}{CG} \\ \cline{3-5}
                       &                             & metis         & k-means        & GNN        &                     \\ \hline
          & 1          & 21    & 9       & 9   & 114  \\
triangles & 2          & 49    & 10      & 10  & 355  \\
          & 3          & 63    & 17      & 17  & 1101 \\ \hline
          & 1          & 46    & 41      & 32  & 655  \\
random    & 2          & 89    & 76      & 53  & 4425 \\
          & 3          & 115   & 102     & 64  & 9893 \\ \hline
          & 1          & 19    & 16      & 15  & 348  \\
Voronoi   & 2          & 36    & 25      & 25  & 987  \\
          & 3          & 39    & 33      & 27  & 3235 \\ \hline
          & 1          & 20    & 17      & 17  & 109  \\
squares   & 2          & 45    & 34      & 34  & 365  \\
          & 3          & 46    & 52      & 52  & 703  \\ \hline
\end{tabular}
\caption[Multigrid iteration counts with different polynomial degrees.]{iteration counts of the MG algorithm to reduce the (relative) residual below $10^{-6}$ employing different initial grids (triangles, random, Voronoi, squares) agglomerated with different strategies (METIS, k-means, GNN) with $\ell = 3$, $p = 1,2,3$, $m = 3$. As a comparison, the iteration counts of the Conjugate Gradient (CG) method are reported in the last column.}
\label{table:MG p}
\end{table}
As expected, since $m$ is fixed, when $p$ increases also the iteration counts increase, because more degrees of freedom are involved in the formulation of the problem \cite{multigridV}.\\
In general, the considered experiments highlight that when employing ML-based agglomeration strategies lower iteration counts are required by the MG solver with respect to METIS, thanks to the higher quality of the produced grids, with GNN having the best performance.

%% file: sections/conclusions.tex
\section{Conclusions}
We presented a new ML-based framework to efficiently handle the open problem of mesh agglomeration for polygonal grids. It is based on the concept of learning the geometrical information contained in the shape of mesh elements, in an automatic and cost-effective manner in order to tailor the approach to a wide range of possible situations, which would not be possible using classical strategies. The novelty of the proposed method consists in improving classical methods for mesh agglomeration using GNN architectures and the k-means clustering algorithm, with GNNs featuring the best performance in terms of accuracy and inference speed. Using these techniques, we can preserve the geometrical quality of the initial grid, enhancing the performance of numerical methods, such as MG solvers, and reducing the mesh complexity in terms of memory storage and associated computational costs. Advantages also include the independence from the differential model and the numerical method at hand, the possibility of handling a wider spectrum of geometries, and the flexibility in exploiting any sort of additional information.\\
\newline
We first re-framed the problem of mesh agglomeration as a graph partitioning problem, by exploiting the graph representation of the connectivity structure of mesh elements. We then developed in this context a framework to employ ML-based strategies, such as k-means and GNNs which can exploit the geometrical information of the grid. We trained GNNs to perform graph partitioning over a suitably constructed database of meshes. We compared the performance of the ML-based strategies with METIS, a classical algorithm for graph partitioning, over a set of different grids, featuring also complex real domains such as a brain MRI scan. Results show that ML-based strategies can better preserve mesh quality, making them suitable for adaptive mesh coarsening. We also employed the proposed algorithms in the context MG methods in PolyDG framework, showing lower iteration counts for ML-based strategies, with GNN having the best performance. Indeed, GNNs have the advantage to process naturally and simultaneously both the graph structure of mesh and the geometrical information, such as the elements areas or their barycentric coordinates. This is not the case with METIS \cite{metis}, which is meant to process only the graph information, or k-means, which can process only the geometrical information. Moreover, GNNs have a significantly lower online computational cost. Being able to compute the solution in a fast manner, even with a slight loss of accuracy, can be particularly beneficial when using multigrid schemes as preconditioners.\\
\newline
Future developments certainly include fostering the generalization capabilities of the GNNs, by generating a larger database of grids or by considering different geometrical features, such as quality metrics or error estimators. Another possibility is to reframe the approach in a reinforcement learning framework, as proposed by \cite{gatti2}. The extension to three-dimensional agglomeration strategies should certainly be explored, where the low online computational cost of GNNs will play a key role. It is worth noticing that all of the proposed agglomeration strategies can be applied to the three-dimensional case with little or no modification, as they mainly rely on the connectivity graph structure of the mesh, which is a dimensionless entity, or on geometrical quantities with natural 
extensions such the area/volume of polytopes of their barycentric coordinates. From a point of view of numerical applications several options can be explored: adaptive mesh coarsening, domain decomposition, using MG schemes as preconditioners, and considering discretization frameworks other than the PolyDG ones such as the one of Virtual Element Methods. Finally, it would be interesting to consider more complex scenarios, for example, geophysical numerical simulations, including fluid-structure interaction with complex and moving geometries. Preliminary results of agglomeration of the brain obtained from medical images are encouraging.

%% file: elsarticle-main.bbl
\begin{thebibliography}{10}
\expandafter\ifx\csname url\endcsname\relax
  \def\url#1{\texttt{#1}}\fi
\expandafter\ifx\csname urlprefix\endcsname\relax\def\urlprefix{URL }\fi
\expandafter\ifx\csname href\endcsname\relax
  \def\href#1#2{#2} \def\path#1{#1}\fi

\bibitem{hyman1997numerical}
J.~Hyman, M.~Shashkov, S.~Steinberg, The numerical solution of diffusion
  problems in strongly heterogeneous non-isotropic materials, Journal of
  Computational Physics 132~(1) (1997) 130--148.

\bibitem{brezzi2005family}
F.~Brezzi, K.~Lipnikov, V.~Simoncini, A family of mimetic finite difference
  methods on polygonal and polyhedral meshes, Mathematical Models and Methods
  in Applied Sciences 15~(10) (2005) 1533--1551.

\bibitem{brezzi2005convergence}
F.~Brezzi, K.~Lipnikov, M.~Shashkov, Convergence of the mimetic finite
  difference method for diffusion problems on polyhedral meshes, SIAM Journal
  on Numerical Analysis 43~(5) (2005) 1872--1896.

\bibitem{da2014mimetic}
L.~Beirao~da Veiga, K.~Lipnikov, G.~Manzini, The mimetic finite difference
  method for elliptic problems, Vol.~11, Springer, 2014.

\bibitem{cockburn2008superconvergent}
B.~Cockburn, B.~Dong, J.~Guzm{\'a}n, A superconvergent ldg-hybridizable
  galerkin method for second-order elliptic problems, Mathematics of
  Computation 77~(264) (2008) 1887--1916.

\bibitem{cockburn2009superconvergent}
B.~Cockburn, J.~Guzm{\'a}n, H.~Wang, Superconvergent discontinuous galerkin
  methods for second-order elliptic problems, Mathematics of Computation
  78~(265) (2009) 1--24.

\bibitem{cockburn2009unified}
B.~Cockburn, J.~Gopalakrishnan, R.~Lazarov, Unified hybridization of
  discontinuous galerkin, mixed, and continuous galerkin methods for second
  order elliptic problems, SIAM Journal on Numerical Analysis 47~(2) (2009)
  1319--1365.

\bibitem{cockburn2010projection}
B.~Cockburn, J.~Gopalakrishnan, F.-J. Sayas, A projection-based error analysis
  of hdg methods, Mathematics of Computation 79~(271) (2010) 1351--1367.

\bibitem{hesthaven2007nodal}
J.~S. Hesthaven, T.~Warburton, Nodal discontinuous Galerkin methods:
  algorithms, analysis, and applications, Springer Science \& Business Media,
  2007.

\bibitem{bassi2012flexibility}
F.~Bassi, L.~Botti, A.~Colombo, D.~A. Di~Pietro, P.~Tesini, On the flexibility
  of agglomeration based physical space discontinuous galerkin discretizations,
  Journal of Computational Physics 231~(1) (2012) 45--65.

\bibitem{antonietti2013hp}
P.~F. Antonietti, S.~Giani, P.~Houston, hp-version composite discontinuous
  galerkin methods for elliptic problems on complicated domains, SIAM Journal
  on Scientific Computing 35~(3) (2013) A1417--A1439.

\bibitem{cangiani2014hp}
A.~Cangiani, E.~H. Georgoulis, P.~Houston, hp-version discontinuous galerkin
  methods on polygonal and polyhedral meshes, Mathematical Models and Methods
  in Applied Sciences 24~(10) (2014) 2009--2041.

\bibitem{antonietti2016review}
P.~F. Antonietti, A.~Cangiani, J.~Collis, Z.~Dong, E.~H. Georgoulis, S.~Giani,
  P.~Houston, Review of discontinuous galerkin finite element methods for
  partial differential equations on complicated domains, in: Building bridges:
  connections and challenges in modern approaches to numerical partial
  differential equations, Springer, 2016, pp. 281--310.

\bibitem{cangiani2017hp}
A.~Cangiani, Z.~Dong, E.~H. Georgoulis, P.~Houston, hp-Version discontinuous
  Galerkin methods on polygonal and polyhedral meshes, Springer, 2017.

\bibitem{antonietti2021high}
P.~F. Antonietti, C.~Facciol{\`a}, P.~Houston, I.~Mazzieri, G.~Pennesi,
  M.~Verani, High--order discontinuous galerkin methods on polyhedral grids for
  geophysical applications: seismic wave propagation and fractured reservoir
  simulations, Polyhedral Methods in Geosciences (2021) 159--225.

\bibitem{beirao2013basic}
L.~Beir{\~a}o~da Veiga, F.~Brezzi, A.~Cangiani, L.~D. Manzini,
  GianM.and~Marini, A.~Russo, Basic principles of virtual element methods,
  Mathematical Models and Methods in Applied Sciences 23~(01) (2013) 199--214.

\bibitem{beirao2014hitchhiker}
L.~Beir{\~a}o~da Veiga, F.~Brezzi, L.~D. Marini, A.~Russo, The hitchhiker's
  guide to the virtual element method, Mathematical models and methods in
  applied sciences 24~(08) (2014) 1541--1573.

\bibitem{beirao2016virtual}
L.~Beir{\~a}o~da Veiga, F.~Brezzi, L.~Marini, A.~Russo, Virtual element method
  for general second-order elliptic problems on polygonal meshes, Mathematical
  Models and Methods in Applied Sciences 26~(04) (2016) 729--750.

\bibitem{da2016mixed}
L.~Beirao~da Veiga, F.~Brezzi, L.~D. Marini, A.~Russo, Mixed virtual element
  methods for general second order elliptic problems on polygonal meshes,
  ESAIM: Mathematical Modelling and Numerical Analysis 50~(3) (2016) 727--747.

\bibitem{beirao2021recent}
L.~Beirao~da Veiga, N.~Bellomo, F.~Brezzi, L.~Marini, Recent results and
  perspectives for virtual element methods, Mathematical Models and Methods in
  Applied Sciences (2021) 1--6.

\bibitem{book.vem.sema.simai.2022}
P.~F. Antonietti, L.~Beir{\~a}o~da Veiga, G.~Manzini, The Virtual Element
  Method and its Applications, Springer International Publishing, 2022.

\bibitem{di2014arbitrary}
D.~A. Di~Pietro, A.~Ern, S.~Lemaire, An arbitrary-order and compact-stencil
  discretization of diffusion on general meshes based on local reconstruction
  operators, Computational Methods in Applied Mathematics 14~(4) (2014)
  461--472.

\bibitem{di2015hybrid}
D.~A. Di~Pietro, A.~Ern, A hybrid high-order locking-free method for linear
  elasticity on general meshes, Computer Methods in Applied Mechanics and
  Engineering 283 (2015) 1--21.

\bibitem{di2015hybrid2}
D.~A. Di~Pietro, A.~Ern, Hybrid high-order methods for variable-diffusion
  problems on general meshes, Comptes Rendus Math{\'e}matique 353~(1) (2015)
  31--34.

\bibitem{di2016review}
D.~A. Di~Pietro, A.~Ern, S.~Lemaire, A review of hybrid high-order methods:
  formulations, computational aspects, comparison with other methods, in:
  Building bridges: connections and challenges in modern approaches to
  numerical partial differential equations, Springer, 2016, pp. 205--236.

\bibitem{di2019hybrid}
D.~A. Di~Pietro, J.~Droniou, The Hybrid High-Order method for polytopal meshes,
  Vol.~19, Springer, 2019.

\bibitem{attene2019benchmark}
M.~Attene, S.~Biasotti, S.~Bertoluzza, D.~Cabiddu, M.~Livesu, G.~Patan{\`e},
  M.~Pennacchio, D.~Prada, M.~Spagnuolo, Benchmark of polygon quality metrics
  for polytopal element methods, arXiv preprint arXiv:1906.01627.

\bibitem{di2021polyhedral}
D.~A. Di~Pietro, L.~Formaggia, R.~Masson, et~al., Polyhedral methods in
  geosciences.

\bibitem{chan1998agglomeration}
T.~F. Chan, J.~Xu, L.~Zikatanov, An agglomeration multigrid method for
  unstructured grids, Contemporary Mathematics 218 (1998) 67--81.

\bibitem{antonietti2020agglomeration}
P.~F. Antonietti, P.~Houston, G.~Pennesi, E.~S{\"u}li, An agglomeration-based
  massively parallel non-overlapping additive schwarz preconditioner for
  high-order discontinuous galerkin methods on polytopic grids, Mathematics of
  Computation.

\bibitem{pan2022agglomeration}
Y.~Pan, P.-O. Persson, Agglomeration-based geometric multigrid solvers for
  compact discontinuous galerkin discretizations on unstructured meshes,
  Journal of Computational Physics 449 (2022) 110775.

\bibitem{gilbert1998geometric}
J.~R. Gilbert, G.~L. Miller, S.-H. Teng, Geometric mesh partitioning:
  Implementation and experiments, SIAM Journal on Scientific Computing 19~(6)
  (1998) 2091--2110.

\bibitem{bassi2012agglomeration}
F.~Bassi, L.~Botti, A.~Colombo, S.~Rebay, Agglomeration based discontinuous
  galerkin discretization of the euler and navier--stokes equations, Computers
  \& fluids 61 (2012) 77--85.

\bibitem{antonietti2015multigrid}
P.~F. Antonietti, M.~Sarti, M.~Verani, Multigrid algorithms for
  hp-discontinuous galerkin discretizations of elliptic problems, SIAM Journal
  on Numerical Analysis 53~(1) (2015) 598--618.

\bibitem{antonietti2017multigrid}
P.~F. Antonietti, P.~Houston, X.~Hu, M.~Sarti, M.and~Verani, Multigrid
  algorithms for hp-version interior penalty discontinuous galerkin methods on
  polygonal and polyhedral meshes, Calcolo 54~(4) (2017) 1169--1198.

\bibitem{xu2017algebraic}
J.~Xu, L.~Zikatanov, Algebraic multigrid methods, Acta Numerica 26 (2017)
  591--721.

\bibitem{chan1998multilevel}
T.~F. Chan, S.~Go, L.~Zikatanov, Multilevel elliptic solvers on unstructured
  grids, in: Computational Fluid Dynamics Review 1998: (In 2 Volumes), World
  Scientific, 1998, pp. 488--511.

\bibitem{multigridV}
P.~Antonietti, G.~Pennesi, V-cycle multigrid algorithms for discontinuous
  galerkinmethods on non-nested polytopic meshes, Journal of Scientific
  Computing 78 (2019) 625--652.

\bibitem{antonietti2007schwarz}
P.~F. Antonietti, B.~Ayuso, Schwarz domain decomposition preconditioners for
  discontinuous galerkin approximations of elliptic problems: non-overlapping
  case, ESAIM: Mathematical Modelling and Numerical Analysis 41~(1) (2007)
  21--54.

\bibitem{antonietti2014domain}
P.~F. Antonietti, S.~Giani, P.~Houston, Domain decomposition preconditioners
  for discontinuous galerkin methods for elliptic problems on complicated
  domains, Journal of Scientific Computing 60~(1) (2014) 203--227.

\bibitem{toselli2004domain}
A.~Toselli, O.~Widlund, Domain decomposition methods-algorithms and theory,
  Vol.~34, Springer Science \& Business Media, 2004.

\bibitem{feng2001two}
X.~Feng, O.~A. Karakashian, Two-level additive schwarz methods for a
  discontinuous galerkin approximation of second order elliptic problems, SIAM
  Journal on Numerical Analysis 39~(4) (2001) 1343--1365.

\bibitem{raissi2019physics}
M.~Raissi, P.~Perdikaris, G.~E. Karniadakis, Physics-informed neural networks:
  A deep learning framework for solving forward and inverse problems involving
  nonlinear partial differential equations, Journal of Computational Physics
  378 (2019) 686--707.

\bibitem{raissi2018hidden}
M.~Raissi, G.~E. Karniadakis, Hidden physics models: Machine learning of
  nonlinear partial differential equations, Journal of Computational Physics
  357 (2018) 125--141.

\bibitem{regazzoni2019machine}
F.~Regazzoni, L.~Ded{\`e}, A.~Quarteroni, Machine learning for fast and
  reliable solution of time-dependent differential equations, Journal of
  Computational Physics 397 (2019) 108852.

\bibitem{regazzoni2020machine}
F.~Regazzoni, L.~Ded{\`e}, A.~Quarteroni, Machine learning of multiscale active
  force generation models for the efficient simulation of cardiac
  electromechanics, Computer Methods in Applied Mechanics and Engineering 370
  (2020) 113268.

\bibitem{hesthaven2018non}
J.~S. Hesthaven, S.~Ubbiali, Non-intrusive reduced order modeling of nonlinear
  problems using neural networks, Journal of Computational Physics 363 (2018)
  55--78.

\bibitem{ray2018artificial}
D.~Ray, J.~S. Hesthaven, An artificial neural network as a troubled-cell
  indicator, Journal of Computational Physics 367 (2018) 166--191.

\bibitem{antonietti2021accelerating}
P.~F. Antonietti, M.~Caldana, L.~Dede, Accelerating algebraic multigrid methods
  via artificial neural networks, arXiv preprint arXiv:2111.01629.

\bibitem{regazzoni2021machine}
F.~Regazzoni, M.~Salvador, L.~Ded{\`e}, A.~Quarteroni, A machine learning
  method for real-time numerical simulations of cardiac electromechanics, arXiv
  preprint arXiv:2110.13212.

\bibitem{ANTONIETTI2022110900}
P.~Antonietti, E.~Manuzzi, Refinement of polygonal grids using convolutional
  neural networks with applications to polygonal discontinuous galerkin and
  virtual element methods, Journal of Computational Physics 452 (2022) 110900.

\bibitem{antonietti2022machine}
P.~Antonietti, F.~Dassi, E.~Manuzzi, Machine learning based refinement
  strategies for polyhedral grids with applications to virtual element and
  polyhedral discontinuous galerkin methods, Journal of Computational Physics
  (2022) 111531.

\bibitem{kmeans}
J.~Macqueen, Some methods for classification and analysis of multivariate
  observations (1967).

\bibitem{hartigan1979algorithm}
J.~A. Hartigan, M.~A. Wong, Algorithm as 136: A k-means clustering algorithm,
  Journal of the royal statistical society. series c (applied statistics)
  28~(1) (1979) 100--108.

\bibitem{likas2003global}
A.~Likas, N.~Vlassis, J.~J. Verbeek, The global k-means clustering algorithm,
  Pattern recognition 36~(2) (2003) 451--461.

\bibitem{bello2012adaptive}
G.~Bello-Orgaz, H.~D. Men{\'e}ndez, D.~Camacho, Adaptive k-means algorithm for
  overlapped graph clustering, International journal of neural systems 22~(05)
  (2012) 1250018.

\bibitem{1706.02216}
W.~L. Hamilton, R.~Ying, J.~Leskovec, Inductive representation learning on
  large graphs (2017).
\newblock \href {http://arxiv.org/abs/arXiv:1706.02216}
  {\path{arXiv:arXiv:1706.02216}}.

\bibitem{gatti1}
A.~Gatti, Z.~Hu, T.~Smidt, E.~G. Ng, P.~Ghysels, Deep learning and spectral
  embedding for graph partitioning (2021).
\newblock \href {http://arxiv.org/abs/arXiv:2110.08614}
  {\path{arXiv:arXiv:2110.08614}}.

\bibitem{gatti2}
A.~Gatti, Z.~Hu, T.~Smidt, E.~G. Ng, P.~Ghysels, Graph partitioning and sparse
  matrix ordering using reinforcement learning and graph neural networks
  (2021).
\newblock \href {http://arxiv.org/abs/arXiv:2104.03546}
  {\path{arXiv:arXiv:2104.03546}}.

\bibitem{xu2021graph}
H.~Xu, Z.~Duan, Y.~Wang, J.~Feng, R.~Chen, Q.~Zhang, Z.~Xu, Graph partitioning
  and graph neural network based hierarchical graph matching for graph
  similarity computation, Neurocomputing 439 (2021) 348--362.

\bibitem{lecun2015deep}
Y.~LeCun, Y.~Bengio, G.~Hinton, Deep learning, nature 521~(7553) (2015)
  436--444.

\bibitem{metis}
G.~Karypis, V.~Kumar, Kumar, v.: A fast and high quality multilevel scheme for
  partitioning irregular graphs. siam journal on scientific computing 20(1),
  359-392, Siam Journal on Scientific Computing 20.
\newblock \href {http://dx.doi.org/10.1137/S1064827595287997}
  {\path{doi:10.1137/S1064827595287997}}.

\bibitem{hu2021posteriori}
X.~Hu, K.~Wu, L.~T. Zikatanov, A posteriori error estimates for multilevel
  methods for graph laplacians, SIAM Journal on Scientific Computing 43~(5)
  (2021) S727--S742.

\bibitem{tarjan1972depth}
R.~Tarjan, Depth-first search and linear graph algorithms, SIAM journal on
  computing 1~(2) (1972) 146--160.

\bibitem{arthur2007proceedings}
D.~Arthur, S.~Vassilvitskii, H.~Gabow, Proceedings of the eighteenth annual
  acm-siam symposium on discrete algorithms, Society for Industrial and Applied
  Mathematics.

\bibitem{urschel2016maximal}
J.~C. Urschel, L.~T. Zikatanov, On the maximal error of spectral approximation
  of graph bisection, Linear and Multilinear Algebra 64~(10) (2016) 1972--1979.

\bibitem{kingma2015adam}
D.~P. Kingma, J.~Ba, Adam: A method for stochastic optimization, Proceedings of
  the 3rd International Conference on Learning Representations (ICLR) (2015)
  1--15.

\bibitem{gary1979computers}
M.~R. Gary, D.~S. Johnson, Computers and intractability: A guide to the theory
  of np-completeness (1979).

\bibitem{andreev2004balanced}
K.~Andreev, H.~R{\"a}cke, Balanced graph partitioning, in: Proceedings of the
  sixteenth annual ACM symposium on Parallelism in algorithms and
  architectures, 2004, pp. 120--124.

\bibitem{karypis1997parmetis}
G.~Karypis, K.~Schloegel, V.~Kumar, Parmetis: Parallel graph partitioning and
  sparse matrix ordering library.

\bibitem{kantabutra2000parallel}
S.~Kantabutra, A.~L. Couch, Parallel k-means clustering algorithm on nows,
  NECTEC Technical journal 1~(6) (2000) 243--247.

\bibitem{ma2018towards}
L.~Ma, Z.~Yang, Y.~Miao, J.~Xue, M.~Wu, L.~Zhou, Y.~Dai, Towards efficient
  large-scale graph neural network computing, arXiv preprint arXiv:1810.08403.

\bibitem{petersen2020neural}
P.~C. Petersen, Neural network theory, University of Vienna.

\end{thebibliography}
